\documentclass[12pt]{article}

  \usepackage{amssymb}         
  \usepackage{amsmath}          
  \usepackage{amsfonts}           
  \usepackage{amsthm}
  \usepackage{amscd}
  \usepackage{enumerate}
  \usepackage[active]{srcltx}
  \usepackage{hyperref}

  \usepackage[mathscr]{eucal}

  \usepackage{graphicx,epsfig}

  \newcommand*{\al}{\alpha}
  
  \newcommand{\Area}{{\operatorname{Area}}}

  \newcommand*{\be}{\beta}

  \newcommand*{\C}{\mathbb{C}}
  \newcommand*{\Chat}{\ensuremath{{\widehat{\mathbb  C}}}}

  \renewcommand*{\d}{{\operatorname{d}}}
  
  \newcommand*{\dA}{\ensuremath{{\operatorname{dA}}}}

  \newcommand*{\D}{\mathbb{D}}

  \newcommand*{\la}{\lambda}

  \newcommand*{\mapfromto}[3]{\hbox{\ensuremath{#1 : #2 \longrightarrow #3}}}

  \newcommand*{\si}{\sigma}

  \newcommand*{\Sm}{\ensuremath{{\setminus}}}

\newtheorem{definition}{Definition}
\newtheorem{theorem}{Theorem}
\newtheorem{corollary}{Corollary}
\newtheorem{lemma}{Lemma}
\newtheorem{proposition}{Proposition}

\newtheorem{remark}{Remark}

\newcommand{\ALIGN}{\begin{align*}}
\newcommand{\ENDALIGN}{\end{align*}}
\newcommand{\ENUM}{\begin{enumerate}}
\newcommand{\ENUMa}{\begin{enumerate}[a.]}
\newcommand{\ENUMA}{\begin{enumerate}[A.]}
\newcommand{\ENUMAF}{\begin{enumerate}[\bf A.]}
\newcommand{\ENUMi}{\begin{enumerate}[i)]}
\newcommand{\ENDENUM}{\end{enumerate}}
\newcommand{\ITMZ}{\begin{itemize}}
\newcommand{\ENDITMZ}{\end{itemize}}
\newcommand{\REFEQN}[1] { \begin{equation}\label{#1} }
\newcommand{\ENDEQN}{\end{equation}}
\newcommand{\THM}{\begin{theorem}}
\newcommand{\REFEXA}[1] { \begin{example}\label{#1} }
\newcommand{\ENDEXA}{\end{example}}
\newcommand{\REM}{ \begin{remark}}
\newcommand{\ENDREM}{\end{remark}}
\newcommand{\REFTHM}[1] { \begin{theorem}\label{#1} }
\newcommand{\RREFTHM}[2] { \begin{theorem}[#1]\label{#2} }
\newcommand{\ENDTHM}{\end{theorem}}
\newcommand{\REFNTH}[1] { \begin{newthm}\label{#1} }
\newcommand{\ENDNTH}{\end{newthm}}
\newcommand{\REFPROP}[1]{\begin{proposition}\label{#1} }
\newcommand{\RREFPROP}[2]{\begin{proposition}[#1]\label{#2} }
\newcommand{\PROP}{\begin{proposition}}
\newcommand{\ENDPROP}{\end{proposition} }
\newcommand{\REFDEF}[1]{\begin{definition}\label{#1} }
\newcommand{\DEF}{\begin{definition}}
\newcommand{\ENDDEF}{\end{definition} }
\newcommand{\REFLEM}[1]{\begin{lemma}\label{#1} }
\newcommand{\RREFLEM}[2]{\begin{lemma}[#1]\label{#2} }
\newcommand{\LEM}{\begin{lemma}}
\newcommand{\ENDLEM}{\end{lemma} }
\newcommand{\REFCOR}[1]{\begin{corollary}\label{#1} }
\newcommand{\RREFCOR}[2]{\begin{corollary}[#1]\label{#2} }
\newcommand{\RCOR}[1]{\begin{corollary}[#1] }
\newcommand{\COR}{\begin{corollary}}
\newcommand{\ENDCOR}{\end{corollary} }
\newcommand{\REFDEFTHM}[1] { \begin{defthm}\label{#1} }
\newcommand{\ENDDEFTHM}{\end{defthm}}

\newcommand{\lemref}[1]{Lemma~\ref{#1}}
\newcommand{\thmref}[1]{Theorem~\ref{#1}}

\newcommand{\PROOF}{\begin{proof}}
\newcommand{\ENDPROOF}{\end{proof}}

\usepackage{color} 
\title{The Intersection of nested filled-in little Julia sets have zero area } 
\author{Carsten Lunde Petersen}

\begin{document}

\maketitle

\begin{abstract}
It is shown that for any infinitely renormalizable quadratic polynomial $Q_c$ 
the intersection $K_\infty$ of all the little filled-in Julia sets (renormalization Julia sets) containing the critical value have zero area.
\end{abstract}


\section{Introduction}\label{intro}
In his lecture during the conference The circle at infinity: an international colloquium in honour of Curtis T. McMullen, Dennis Sullivan gave a lecture in which all drawings were turned in unusal directions, e.g. the Mandelbrot set was drawn with the tail pointing down and so forth. It created a certain intended frustration for some, but Sullivans point was clear:
After nearly fifty years work on and around the Mandelbrot set, a new break through seems to need completely new ideas, perhaps coming from minds, which do not carry the burden of 50 years intense developments, and so are able to look at matters with fresh eyes. 
One potential such new idea was published September 20.~this year on ArXiv by Zihao Ye, 
who in dialog with ChatGPT had found a completely new and elementary proof of Sullivans celebrated no-wandering domains theorem for rational maps proved in the 1980-ties using Teichm{\"u}ller theory, 
\cite{Sullivan-nwdt}. 
Indeed Ye's result has already been used by Prochorov, Rempe and Waterman, likewise in collaboration with generative AI, to solve an 80 years old conjecture, on absence of bounded-orbit wandering domains for entire transcendental maps, \cite{Prochorov-Rempe-Waterman} . 
Ye's proof roughly one page long exploits in a clever way that the Riemann sphere with finitely many (but at lest $3$) points removed has finite hyperbolic area given by the number of removed points. 
Perhaps the most amazing property of the proof is that it has gone un-noticed for so long, underlining Sullivans point above. It only uses mathematical entities, which where 
available already to Fatou and Julia.

This short note created purely by use of HI (human inteligence) shows that exactly the same idea 
gives that, if a quadratic polynomial is infinitely renormalizable, then the nested intersection $K_\infty$ of the renormalization Julia sets of all the renormalizations around the critical value has zero Euclidean Area. Clearly $K_\infty$ is connected, being the intersection of nested connected sets. 

There are many cases in which all the renormalizations have a robust fundamental annulus 
and so $K_\infty$ is reduced to a point. 
But there are also examples where $K_\infty$ is non locally connected, and in fact with diameter arbitrarily close to $1$, \cite{Sorensen}.
And there are examples where the full Julia sets have positive Area, 
\cite{Buff-Cheritat}, \cite{Avila-Lyubich}. 
So the take home message here is that such positive Area must come from the set of points, 
whose orbits accumulates $K_\infty$ without ever hitting it. 
Of course the main question still remains, 
are there infinitely renormalizable quadratic polynomials, whose Julia set carries 
a measureable line field with positive Area support. 
The conjecture MLC, the Mandelbrot set is locally connected, 
if true implies that such strange quadratic polynomials do not exist 
and equivalently every connected component of the interior of the Mandelbrot set is hyperbolic. 

To keep this note short I proceed straight to the point referring the uninitiated reader to the 
substantial literature on holomorphic dynamics and the Mandelbrot set, 
\cite{Milnor_Dynamics, Carleson-Gamelin, Steinmetz, orsaynotesI, orsaynotesII, Milnorasterisque, polylikemaps, McMBook1, universalityofM, Milnorquad} m.m.

The notion of polynomial-like maps was introduced by Douady and Hubbard in the landmark paper 
\cite{polylikemaps}. A polynomial-like map is a proper holomorphic map {\mapfromto f U {U'}}
of degree $d\geq 2$, where $U\subset\subset U'$ are Riemann surfaces isomorphic to $\D$. 
When $d=2$ and $U'\subset\C$ we shall call such a map a \emph{quadratic-like} map. 
The \emph{filled-in Julia set} of $f$ is the compact set $K_f := \{z\in U \mid \forall n\ge 0 : f^n(z)\in U\}$. 
A hybrid equivalence between two polynomial-like maps is a quasi-conformal homeomorphism, 
which conjugates dynamics of the maps possibly reducing them to polynomial-like restrictions preserving the filled-in Julia sets, and which is conformal $a.e.$ between the filled-in Julia sets. 
According to the Straightening Theorem \cite[Theorem 1]{polylikemaps}, 
any polynomial-like map of degree $d$ is hybridly- equivalent to the restriction of an actual polynomial of degree $d$, which is unique up to affine conjugacy, when $K_f$ is connected. 
A quadratic polynomial has counting multiplicity $2$ fixed points, one of which called 
$\be$ is the landing point of the external ray of argument $0$, the other fixed point is called $\al$.
For a quadratic-like map $f$ with $K_f$ connected, 
the straightening theorem uniquely induces labels $\be_f$ and $\al_f$ of the fixed points of $f$.

A quadratic-like map {\mapfromto f U {U'}} is called renormalizable if there exists $n>1$ and a quadratic-like restriction {\mapfromto {F=f^n} V {V'}} 
with connected little filled-in Julia set $K_F$ and $f^j(V)\subset U$ for $0\le j < n$. 
We can and shall suppose the critical value of $f$ belongs to $K_F$ and so equals the critical value of $F$. Moreover at the expense of replacing $n$ by $qn$ for some natural number $q$, 
we can and shall assume the little filled-in Julia sets 
$K_F, f(K_F), \ldots f^{n-1}(K_F)$ are disjoint except possibly for $\be_F$, 
which can be a common point, in which case we call the renormalization a \emph{satellite renormalization}. Otherwise it is called a primitive renormalization. 
We shall say that $F$ is a top-level renormalization if $n$ is minimal. 
If $F$ is a top-level satellite renormalization then $\be_F = \al_f$. 
In all cases satellite or primitive, the fixed point $\be_f$ is strongly separated from $K_F$. 
It follows that if $F$ is also renormalizable with connected little filled-in Julia set $K\subset K_F$, 
then the first $n$ sets $K, f(K), \ldots f^{n-1}(K)$ are disjoint and $f^j_{|K}$ is a holomorphic 
diffeomorphism onto the image for $0\le j <n$.

The polynomial $Q_c(z) = z^2+c$, $c\in\C$ is called infinitely renormalizable, 
if there exists a sequence $\{n_k\}_k$ $n_k>1$ of relative periods and renormalizations 
\REFEQN{infinte_renormalizations}
{\mapfromto {f_k=Q_c^{N_k}} {U_k} {U_k'}},\qquad\textrm{where}\qquad N_k = \prod_{j=1}^k n_j
\ENDEQN
and $U_k\subset U_{k+1}$, $U_k'\subset U_{k+1}'$. Set $K_\infty := \cap_{k=1}^\infty K_{f_k}$ 
and denote by $\Area(\cdot)$ the Euclidean area measure. Then

\REFTHM{main}
Suppose $Q_c$ is infinitely renormalizable as defined above then 
$$
\Area(K_\infty) = 0.
$$
\ENDTHM

We shall apply the method of Ye, \cite{Ye} in the interpretation of Buff \cite{Buff}. 

Fix an infinitely renormalizable quadratic polynomial $Q_c$ with renormalizations as in 
\eqref{infinte_renormalizations}.

Let $X_n$ denote the set of periodic points for $Q_c$ with periods at most $n$. 
These are all repelling and $K_\infty\cap X_n=\emptyset$. 
Set $Y_n := X_n\cup\{c\}$ and $Z_n := Q_c^{-1}(Y_n)$, so that $X_n\subset Y_n\cap Z_n$ and 
the restriction {\mapfromto {Q_c} {Z_n}{Y_n}} is a covering map.
Then $|X_n| = 2^n$, so that the complements in $\C$ of all three sets are hyperbolic. 
Let $\tau_n, \si_n$ and $\rho_n$ denote the hyperbolic metrics with curvature $-1$ 
on $\C\Sm Z_n , \C\Sm Y_n \subset \C\Sm X_n$ respectively. 
For $\la\in\{\tau_n, \si_n, \rho_n\}$ denote also by $\la$ the coefficient function of that 
hyperbolic metric relative to the Euclidean metric and write $\dA_\la = \la(z)^2\dA$ for 
the corresponding area elements, where $\dA$ is the Euclidean area element. 
Then by Gauss-Bonnet's theorem
$$
\int_{\C\Sm X}\dA_\la = 2\pi(|X|-1),
$$ 
where $|X|$ denotes the number of elements in $X$. 
Moreover $\dA_{\rho_n}< \dA_{\si_n}$ on $\C\Sm Y_n$ and
$\dA_{\rho_n} < \dA_{\tau_n} = \dA_{Q_c^*\si_n}$ on $\C\Sm Z_n$.

The measure $\nu_n$ given by $\d\nu_n := \dA_{\si_n} - \dA_{\rho_n} = (\si_n^2-\rho_n^2)\dA_s$ is a non-negative measure on $\C\Sm Y_n$ (on $\C$) with 
\REFEQN{absolute_bound}
\nu_n(\C\Sm Y_n) = \int_{\C\Sm Y_n} (\si_n^2-\rho_n^2)\dA_s = 2\pi(|Y_n|-1 - (|X_n|-1)) = 2\pi.
\ENDEQN
Set $K' = \cup_{j=1}^\infty Q_c^j(K_\infty)$. 
Then $K_\infty\cap K' = \emptyset$ and $Q_c$ maps $K_\infty\cup K'$ diffeomorphically onto $K'$. 
That is $K_\infty$ is a wandering compact subset avoiding all periodic points as well as the unique critical point $0$ More precisely for every $j\geq 1$ the restriction {\mapfromto {Q_c} {Q_c^j(K_\infty)}{Q_c^{j+1}(K_\infty)}} is a diffeomorphism and for $0\le m < n$ :
$Q_c^m(K_\infty)\cap Q_c^n(K_\infty)=\emptyset$.  
\REFLEM{hyperbolic_area_lemma}
$$
\int_{K_\infty}\dA_{\rho_n} \le  \int_{K'}\d\nu_n \le 2\pi.
$$
\ENDLEM
See \cite[Lemme 1]{Buff} and \cite[Lemma 3.2]{Ye} for the corresponding statement and proof in the case of a hypothetical wandering component of a rational map.
\PROOF
\ALIGN
\int_{K_\infty}\dA_{\rho_n} + \int_{K'}\dA_{\rho_n} = \int_{K_\infty\cup K'}\dA_{\rho_n}
&\le  \int_{K_\infty\cup K'}\dA_{Q_c^*\si_n}\\
&= \int_{K'}\dA_{\si_n} = \int_{K'}\d\nu_n + \int_{K'}\dA_{\rho_n}
\end{align*}
Since $\int_{K'}\dA_{\rho_n} \le \int_{\C\Sm X_n}\dA_{\rho_n} = 2^n\cdot 2\pi < \infty$, 
we can cancel this integral from both sides of the inequality and the Lemma then follows from 
\eqref{absolute_bound}.
\ENDPROOF
\thmref{main} now follows. Indeed $K_\infty  = \cap_{k=1}^\infty K_{f_k} \subset J_c$ and repelling periodic points are dense in $J_c$ so that  
$$
\inf\{\rho_n(z)\mid z\in J_c\supset K_\infty\} \to\infty\qquad\textrm{as}\qquad n\to\infty.
$$ 
Hence if $\Area(K_\infty)>0$, then there would exist $n$ such that 
$\int_{K_\infty}\dA_{\rho_n} > 2\pi$, which contradicts that these hyperbolic areas are all less than $2\pi$ by \lemref{hyperbolic_area_lemma}.

\REM
Let $X\subset\Chat$ be a finite set with at least $3$ distinct points, so that $\Chat\Sm X$ is 
hyperbolic. Let $\la$ denote both the hyperbolic metric and the coefficient function of the hyperbolic metric relative to the spherical metric, which in turn on $\C$ has the coefficient function 
$s(z) = 2/(1+|z|^2)$ relative to the Euclidean metric.
Then the hyperbolic area measure $A_\la$ is absolutely continuous with respect 
to the spherical area measure $\Area_s$.
Indeed $z\mapsto 1/z$ is a spherical isometry and $\la^2$ is integrable with respect to $\Area_s$ 
so that 
$$
\int_\Chat \dA_\la = \int_{\C\Sm X}\dA_\la = 2\pi(|X|-2),
$$ 
where $|X|$ denotes the number of elements in $X$. 
\ENDREM

Address:

Carsten Lunde Petersen, 
Department of Mathematical Sciences, 
University of Copenhagen, 
Universitetsparken 5, 
DK-2100 Copenhagen {\O}.
e-mail: lunde@math.ku.dk


\begin{thebibliography}{MMM}
\bibitem[A-L]{Avila-Lyubich}
A.~Avila and M.~Lyubich ---
\emph{Lebesgue measure of Feigenbaum Julia sets}, 
Ann.~of Math. {\bf 195} No. 1 (2022), pp.1-88, 
DOI: 10.4007/annals.2022.195.1.1

\bibitem[B-C]{Buff-Cheritat} ---
X.~Buff and A.~Cheritat, \emph{Quadratic Julia sets with positive area}, 
Ann.~of Math. {\bf 176} No. 2 (2012), pp. 673–746,
https://doi.org/10.4007/annals.2012.176.2.1.

 \bibitem[B]{Buff}{X.~Buff}
 X.~Buff, --- \emph{ Th{\'e}or{\`e}me de non errance},
Lecture notes, September 2026.

\bibitem[C-G]{Carleson-Gamelin} 
L.~Carleson and T.~Gamelin, \emph{Complex Dynamics}, Springer-Verlag 1993.

 \bibitem[DH1]{orsaynotesI} A. Douady, J.H. Hubbard --- 
\emph{Etude dynamique des polyn\^omes complexes, I}, Publications Math\'ematiques
d'Orsay, 84-02, 1984.

\bibitem[DH2]{orsaynotesII} A. Douady, J.H. Hubbard --- 
\emph{Etude dynamique des polyn\^omes complexes, II}, publications math\'ematiques
d'Orsay, 85-01, 1985.

\bibitem[DH3]{polylikemaps} A. Douady, J.H. Hubbard --- 
\emph{On the dynamics of polynomial-like mappings},
Ann. Scient. Ec. Norm. Sup., t.18, p. 287-343, 1985.

\bibitem[McM1]{McMBook1}C.~McMullen --\emph{Complex Dynamics and Renormalization.}
 Annals of Math. Studies, {\bf 135}, Princeton University Press, Princeton, NJ, 1995. 

\bibitem[McM2]{universalityofM}C.~McMullen --\emph{The Mandelbrot set is universal.}
The Mandelbrot set, theme and variations pp 1--17, London Math. Soc. Lecture Note Ser., 
274, Cambridge Univ. Press, Cambridge, 2000.

\bibitem[Mi1]{Milnorquad} J. ~Milnor --- 
\emph{Geometry and Dynamics of quadratic Rational Maps}, Experimental Mathematics {\bf 2}:1, (1993), pp. 37-83.

\bibitem[Mi2]{Milnorasterisque} J.~Milnor --- 
 \emph{Periodic orbits, external rays and the Mandelbrot set : an expository account.}
Asterisque 261 (2000) `Geometrie Complexe et Systemes Dynamiques', pp. 277-333.

\bibitem[Mi3]{Milnor_Dynamics} 
J.~Milnor, \emph{Dynamics in One Complex Variable}. 
Annals of Mathematics Studies, vol. 160, (2006)Princeton University Press, 

\bibitem[P-R-W]{Prochorov-Rempe-Waterman} 
N.~Prochorov, L.~Rempe and J.~Waterman ---
\emph{Absence of bounded-orbit wandering domains},
https://doi.org/10.48550/arXiv.2609.28279

\bibitem[St]{Steinmetz} N.~Steinmetz, 
\emph{Rational Iteration: Complex Analytic Dynamical Systems}, 
Berlin, New York: De Gruyter, 1993. https://doi.org/10.1515/9783110889314

\bibitem[Su]{Sullivan-nwdt} D.~Sullivan, 
\emph{Quasiconformal homeomorphisms and dynamics I. Solution of the Fatou–Julia problem on
wandering domains}, Ann.~of Math. {\bf 122} No. 2 (1985), , pp. 401–418, 
https://www.jstor.org/stable/1971308.

\bibitem[S{\o}]{Sorensen} D.~E.~K.~S{\o}rensen,
\emph{Infinitely Renormalizable Quadratic Polynomials, with Non-Locally Connected Julia Set}, 
Geom Anal {\bf 10}, (2000) pp.~169–206, 
DOI: 10.1007/BF02921810

\bibitem[Y]{Ye} {Zihao Ye} --- 
\emph{A new proof of no-wandering domains theorems without quasiconformal techniques}. arXiv:2609.23834  
\end{thebibliography}
\end{document}